\documentclass[a4paper, 12pt,reqno]{amsart}

\usepackage[T1]{fontenc}
\usepackage[utf8]{inputenc}
\usepackage[english]{babel}
\usepackage{url} 
\usepackage{amsmath}
\usepackage{amssymb}
\usepackage{amsthm}
\usepackage{mathtools}
\usepackage[italian]{varioref}
\usepackage{hyperref}
\usepackage{typehtml}
\usepackage{multicol}
\usepackage{color}
\usepackage{mathrsfs}
\usepackage{comment}
\usepackage{mathabx}
\usepackage{tikz-cd}
\usepackage{adjustbox}
\usepackage{bbm}
\usepackage{hyperref}
\usepackage[a4paper,top=3cm,bottom=3cm,left=3.5cm,right=3.5cm]{geometry}

\usepackage{scalerel,stackengine}
\stackMath
\newcommand\reallywidehat[1]{%
\savestack{\tmpbox}{\stretchto{%
  \scaleto{%
    \scalerel*[\widthof{\ensuremath{#1}}]{\kern-.3pt\bigwedge\kern-.3pt}%
    {\rule[-\textheight/2]{2ex}{\textheight}}
  }{\textheight}%
}{1ex}}%
\stackon[1pt]{#1}{\tmpbox}%
}
\hypersetup{hidelinks}

\newcommand{\R}{\mathbb R}

\newcommand{\T}{\mathbb T}
           
\newcommand{\e}{\mathrm{e}}                
\newcommand{\de}{\mathop{}\!\mathrm{d}}

\newcommand{\pa}{\partial}

\newcommand{\fcr}{\mathcal X}

\def\ie{\textit{i.e.} }
\def\be{\begin{equation}}
\def\ee{\end{equation}}
\def\bea{\begin{eqnarray}}
\def\eea{\end{eqnarray}}

\def\nn{\nonumber}

\makeatletter
\@namedef{subjclassname@2020}{\textup{MSC2020} subject classifications}
\makeatother

\newtheorem{thm}{Theorem}

\numberwithin{equation}{section}

\title[Propagation of chaos for topological interactions]
{Propagation of chaos for topological interactions by a coupling technique}

\author[P. \ Degond]{Pierre Degond}
\address{Pierre Degond \hfill\break \indent
	Institut de Math\'ematiques de Toulouse ; UMR5219
	\hfill\break \indent
    	Universit\'e de Toulouse ; CNRS \hfill\break \indent
	UPS, F-31062 Toulouse Cedex 9, France}
\email{pierre.degond@math.univ-toulouse.fr}

\author[M. \ Pulvirenti]{Mario Pulvirenti}
\address{Mario Pulvirenti\hfill\break \indent
	Dipartimento di Matematica, 
	\hfill\break \indent
	Universit\'a di Roma `La Sapienza'
          \hfill\break \indent
          Piazzale Aldo Moro 2, 00185 Roma, Italy}
\address{International Research Center M\&MoCS, 
	\hfill\break \indent
	Universit\'a degli studi dell’Aquila
          \hfill\break \indent
          Via Giovanni Gronchi 18, 67100 L'Aquila, Italy}
\address{Accademia Nazionale dei Lincei
	 \hfill\break \indent   
	Via della Lungara 10,
	00165 Roma, Italy}
\email{pulviren@mat.uniroma1.it}

\author[S.\ Rossi]{Stefano Rossi}
\address{Stefano Rossi\hfill\break \indent
	Dipartimento di Matematica,
          \hfill\break \indent
	 Universit\`a di Roma `La Sapienza'
	\hfill\break \indent
	Piazzale Aldo Moro 2, 00185 Roma, Italy}
\email{stef.rossi@uniroma1.it}

\begin{document}

\date{\today}
\begin{abstract}
We consider a system of particles which interact through a jump process. The jump intensities are functions of the proximity rank of the particles, a type of interaction referred to as topological in the literature. Such interactions have been shown relevant for the modelling of bird flocks. We show that, in the large number of particles limit and under minimal smoothness assumptions on the data, the model converges to a kinetic equation which was derived in earlier works both formally and rigorously under more stringent regularity assumptions. The proof relies on the coupling method which assigns to the particle and limiting processes a joint process posed on the cartesian product of the two configuration spaces of the former processes. By appropriate estimates in a suitable Wasserstein metric, we show that the distance between the two processes tends to zero as the number of particles tends to infinity, with an error typical of the law of large numbers. 
\end{abstract}


\maketitle

\pagestyle{headings}

\section{Introduction}

Systems of self-propelled agents undergoing local interactions are ubiquitous in nature, from migrating cells \cite{GBKM} to locust swarms \cite{BBHA} and fish schools \cite{LGCT}. They form intriguing patterns such as coherent motion, travelling bands, oscillations etc. encompassed in the generic term of collective dynamics (see a review in \cite{VZ}). Most models of collective dynamics are based on mean-field interactions (such as the Cucker-Smale \cite{CS} or Vicsek \cite{VCBCS} models) or binary contact interactions \cite{BDG}. However, a third type of interaction has been suggested following observations of bird flocks \cite{BCCC,CCGP} and referred to as ``topological interaction''. In this kind of interaction, the strength of the interaction of an agent with another one is a function of the proximity rank of the latter with respect to the former. The seminal paper \cite{BCCC} has been followed by a number of papers studying various aspects of this phenomenon see e.g. \cite{BFW,CCGPS,GC,NMG, SB}. 

Mathematically, flocking of systems of topologically interacting particles have been investigated in \cite{Ma,SB2,WC}. In \cite{H}, in addition to studying flocking, the author proposes kinetic and fluid models derived from mean-field topological interactions. The present work is strongly aligned with \cite{BD,BD2,DP} where kinetic models are derived for topological interaction models based on jump processes. More precisely, \cite{DP} proves propagation of chaos and provides a rigorous proof of the model formally derived in \cite{BD}. The proof of \cite{DP} makes the limiting assumption that the interaction strength is an analytic function of the normalized rank (a concept precisely defined below) and is based on the BBGKY hierarchy. In the present work, we propose an alternative proof of the result of \cite{DP} based on the coupling method. The advantage of the coupling method over the BBGKY hierarchy is that it only requires the interaction strength to be Lipschitz continuous, a much more general and natural assumption than that of \cite{DP}. On the other hand, \cite{BD2} formally derives a kinetic model for a more singular interaction. The mathematical validity of this formal result is still open. The literature on propagation of chaos and derivation of kinetic models from particle ones is huge and it is difficult to provide a fair account of all relevant contributions in a short introduction. We refer the interested reader to the reviews 
\cite{CD1,CD2} which provide a fairly detailed description of the subject.

The outline of this paper is as follows. In Section \ref{sec:presentation}, we present the model and provide a formal derivation of the macroscopic model. We then state the theorem and comment it in view of the previous results. Section \ref{sec:proof} is devoted to the proof.

\section{Presentation of the model and main results}
\label{sec:presentation}

We recall the model and notations introduced in \cite{BD, DP} 
and state our result. We study a $N$-particle system in $\R^d$, 
$d=1,2,3 \dots$ ( or in $\T^d$ the $d$-dimensional torus).
Each particle, say particle $i$, has a position $x_i$ and velocity $v_i$. 
The configuration of the system is denoted by
$$
Z_N=\{z_i\}_{i=1}^N=\{(x_i,v_i) \}_{i=1}^N= (X_N,V_N).
$$ 

Given the particle $i$, we order the remaining particles 
$j_1,j_2, \cdots j_{N-1}$ according to their distance from $i$, 
namely by the following relation
$$
| x_i -x_{j_h} | \leq | x_i -x_{j_{h+1}} |, \qquad h=1,2 \cdots N-1.
$$

The rank $R(i,k)$ of particle $k=j_h$ (with respect to  $i$) is $h$. 
Note that, 
if $B_r(x)$ denotes the closed ball of center $x \in \R^d$ and radius $r>0$,
we have
$$
R(i,k)=\sum\limits_{\substack {1\leq h \leq  N \\ h\neq i}} \fcr_{B_{|x_i -x_k|}(x_i)}(x_h),
$$
where $\fcr_A$ is the characteristic function of the set $A$.

Given a non-increasing Lipschitz continuous function
$$
K: [0,1]\to \R^+  \quad \text{s.t.} \quad \int_0^1 K(r) \de r =1,
$$
we introduce the transition probabilities
\be
\label{prob}
\pi^N_{i,j} = \frac {K(r(i,j))} {\sum_s K( \frac s {N-1} )},
\ee
where $r(i,j)$ is the normalized rank:
$$
r(i,j)= \frac {R(i,j)}{N-1} \in \Big\{ \frac 1{N-1}, \frac 2{N-1} , \dots \Big\}.
$$

Thanks to the normalization in \eqref{prob}, we have that 
$
\sum_{j}^{} \pi^N_{i,j}=1.
$
We can also rewrite $\pi^N_{i,j}$ as
\begin{equation}
\label{massa_part}
\pi^N_{i,j}= \alpha_N K \Big( r(i,j) \Big),
\end{equation}
where
\be
\label{alpha}
\alpha_N= \frac{1}{(N-1)(1- e_K(N))}
\ee
and $e_K(N)$ is the error given by the Riemann sums
\be
\label{riemann}
e_K(N)=\int_0^1 K(x) \de x - \frac{1}{N-1}\sum_s K\Bigl(\frac{s}{N-1}\Bigr).
\ee

We are now in position to introduce a stochastic process 
describing alignment via a topological interaction.
The particles go freely: $ x_i+v_i t$. 
At some random time dictated by a Poisson process of intensity $N$, 
choose a particle (say $i$) with probability $\frac 1N$ 
and a partner particle, say $j$, with probability $\pi_{i,j}$. 
Then perform the transition
$$
(v_i,v_j) \to (v_j, v_j).
$$
After that the system goes freely with the new velocities and so on. 

The process is described by the following Markov generator given, 
for any $\Phi \in C^1_b (\R^{2dN})$, by
\begin{align}
\label{generator1}
L_N \Phi ( X_N,V_N) &=   \sum\limits_{i=1}^{N} v_i \cdot \nabla_{x_i} \Phi (X_N, V_N) \\
& +\sum\limits _{i=1}^N \sum\limits_{\substack {1\leq j \leq  N \\ i\neq j}} \pi^N_{i,j}
\big[ \Phi (X_N, V^i_N(v_j)) - 
 \Phi (X_N, V_N)  \big], \nn
\end{align}
where $V^i_N(v_j)=(v_1 \dots v_{i-1}, v_j, v_{i+1} \dots v_N)$ if $V_N=(v_1 \dots v_{i-1}, v_i, v_{i+1} \dots v_N)$.

Note that $\pi_{i,j}^N$ depends not only on $N$ but also on the whole spatial configuration $X_N$. Therefore
the law of the process $W^N(t)=W^N(Z_N;t)$ is driven by the following evolution equation
\begin{align}
\label{master}
\pa_t  \int W^N(t) \Phi &= \int W^N (t) \sum\limits_{i=1}^{N} v_i \cdot \nabla_{x_i} \Phi  \\
&
+\int W^N (Z_N;t) 
  \sum\limits _{i=1}^N  
  \sum\limits_{\substack {1\leq j \leq  N \\ i\neq j}} \pi^N_{i,j}\big[ \Phi (X_N, V_N^i(v_j)) 
  -\Phi (X_N, V_N)  \big] ,  \nn
\end{align}
for any test function $\Phi$.

We assume that the initial measure $W^N(0)$ factorizes, namely $W^N(0)= f_0^{\otimes N} $
where $f_0$ is the initial datum for the limiting kinetic equation we are going to establish. 
Note also that $W^N(Z_N;t)$, for $t \geq 0$,  is symmetric in the exchange of particles.

The strong form of equation \eqref{master} is
$$
\Big(\pa_t +\sum\limits_{i=1}^{N} v_i \cdot \nabla_{x_i}\Big)W^N(t) =-N W^N(t)+{\mathcal L}_NW^N(t)
$$
where
$$
{\mathcal L}_N W^N(X_N,V_N;t)= \sum\limits _{i=1}^N  \sum\limits_{\substack {1\leq j \leq  N \\ i\neq j}} \int du \, \pi^N_{i,j} \,
 W^N (X_N, V_N^{(i)}(u);t) \delta (v_i-v_j).
$$

\subsection{Heuristic derivation}

We now want to derive the kinetic equation we expect to be valid in the limit $N \to \infty$.
Setting $\Phi(Z_N)=\varphi (z_1)$ in \eqref {master},  we obtain
\be
\label{kin1}
\pa_t \int f^N_1 \varphi = \int f^N_1 v \cdot \nabla_{x} \varphi -\int f^N_1 \varphi  
+\int W^N \sum_{j\neq 1} \pi^N_{i,j} \varphi (x_1, v_j).
\ee
Here $f^N_1$ denotes the one-particle marginal of the measure $W^N$. We recall that the $s$-particle marginals are defined by
\be
\label{marginals}
f^N_s (Z_s) =\int W^N (Z_s, z_{s+1} \cdots z_N) dz_{s+1} \cdots dz_N, \qquad s=1,2 \cdots N
\ee
and are the distribution of the first $s$ particles (or of any group of $s$ tagged  particles).

In order to describe the system in terms of a single kinetic equation, we expect that chaos propagates.  
Actually since $W^N$ is initially factorizing, although the dynamics creates correlations, we hope that, due to the weakness of the interaction, 
factorization still holds approximately also at any positive time $t$, namely
$$
f^N_s \approx f_1 ^{\otimes s}.
$$
In this case the law of large numbers does hold, that is
$$
\frac 1N \sum_j \delta (z-z_j) \approx f^N_1(z,t) 
$$
for $W^N$- almost all $Z_N=\{z_1 \cdots z_N\} $. Then
\begin{align*}
\pi^N_{i,j} &\approx  \frac 1{N-1} K \Big( \frac 1{N-1}  \sum_k \fcr_{B_{|x_i-x_j|}(x_i)} (x_k) \Big ) \nn \\
&\approx \frac 1{N-1} K \Big( M_{\rho} ( B_{|x_1-x_2|}(x_1)) \Big)
\end{align*}
where
\begin{equation}
\label{mass}
M_{\rho} (B_R(x))= \int_{B_R(x)} \rho(y)\de y,
\end{equation}
and $\rho (x)  =\int \de v f^N_1 (x,v) $ is the spatial density.
Motivated by this remark, from now on we use the following notation
$$
M_{X_N}(B_{|x_i - x_j|}(x_i))= r(i,j)=\frac 1{N-1}  \sum_k \fcr_{B_{|x_i-x_j|}(x_i)} (x_k).
$$
Here  $M$ stands for `mass' and the notation introduced is justified by the law of large numbers.

In conclusion we expect that, by \eqref {kin1}, in the limit $N\to \infty$, $f^N_1 \to f$ and $f^N_2 \to f^{\otimes 2} $, where $f$ solves
\begin{equation*}
\pa_t \int f \varphi = \int f v \cdot \nabla_{x} \varphi -\int f \varphi  
+\int f(z_1) f(z_2)  \varphi (x_1, v_2) K\Big( M_{\rho} ( B_{|x_1-x_2|}(x_1))\Big)
\end{equation*}
which is the weak form of the equation
\be
\label{eqkin}
\Big(\pa_t + v \cdot \nabla_x \Big)f(x,v,t)= -f(x,v,t) + \rho(x,t) \int K \Big(M_{\rho} ( B_{|x-y|}(x) )\Big) \,f(y,v,t) \de y.
\ee
We remark that existence and uniqueness of global solutions in $L^1(\R^{2d})$ for the kinetic equation \eqref{eqkin}
can be proved by using a standard Banach fixed-point argument.

Once known $f$, we can construct the one-particle nonlinear process given by the generator
\[
L^{(1)}_1 \phi(x,v)=(v\cdot \nabla_x -1) \phi(x,v) + \int f(y,w) \phi(x,w) K \Big( M_\rho(B_{|x-y|}(x) ) \Big) \de y \de w.
\]
We also introduce the $N$-particle process given by $N$ independent copies of the above process. 
Its generator is
\begin{multline}
\label{generator2}
L^{(1)}_N \Phi(Z_N)=V_N \cdot \nabla_{X_N} \Phi(Z_N) \\
+ \sum_i \Big [ \int \Phi(X_N, V_N^i(w_i))K\Big(M_\rho( B_{|x_i - y_i|}(x_i)) \Big)f(y_i, w_i)\de y_i \de w_i  - \Phi(X_N, V_N) \Big ].
\end{multline}

\subsection{Motivations and main result}

This work aims to prove propagation of chaos for the $N$-particle process described by \eqref{generator1}.
Propagation of chaos consists in preparing a system of $N$ particles with initial configurations
i.i.d with a given law $f_0$
and show that,
considering any group of fixed $s$ particles between the $N$ ones,
this independence (chaos) is also recovered for future times for the fixed $s$-group when $N \to \infty$. 
This is expressed mathematically by saying that the $s$-particle marginal 
$f^N_s(t)$ introduced in \eqref{marginals}
approximates $f^{\otimes s}(t)$ for positive times, 
where $f(t)$ is the solution with initial datum $f_0$ 
of the limit equation \eqref{eqkin}.

As mentioned in the introduction,
the propagation of chaos result for \eqref{generator1} was already obtained 
in \cite{DP} using hierarchical techniques.
Indeed, the BBGKY hierarchies are a powerful 
approach 
but their structure is such that the equation for the $s$-marginal
depends only on the $(s+1)$-marginal. 
In this case the non-binary nature of the topological interaction 
does not allow to derive this hierarchical structure, 
unless the interaction function $K$ is real analytic and therefore expandable in series,
which is exactly the assumption made in \cite{DP}.

The reason for this work is to provide a different derivation of the limit kinetic equation,
using the classic probabilistic coupling technique.
In general, given two stochastic processes $X$ and $Y$,
a coupling is a realization of a new process
on a product probability space that has as marginal distributions those of $X$ and $Y$.
This approach brings a more natural and general proof,
avoiding the analyticity assumption on $K$.

\begin{thm}

Let $f \in C([0,T]; L^1(\R^{2d}))$ solution of the limit equation \eqref{eqkin} with initial datum $f_0 \in L^1(\R^{2d})$.
Assume that the interaction function $K$ is Lipschitz-continuous 
and consider the $N$-particle dynamics such that $W_N(0)=f_0^{\otimes N}$.

If $f^N_s$ denotes the $s$-marginal as defined in \eqref{marginals}, 
for $t \in [0,T]$ and $s \in \{1, \dots, N\}$,
it holds that
\be
\label{main}
\| f^N_s(t) - f ^{\otimes s}(t)\|_{L^1(\R^{2ds})} \le s\frac{ \e^{C_K T}}{\sqrt{N-1}},
\ee
where $C_K$ is a constant depending only on the Lipschitz constant of $K$.
\end{thm}

The topological character of the interaction
bring us naturally to work with norms of strong type
and in particular with the $L_1$/Total variation distance
(see also \cite{BCR} 
where a distance similar to the Total Variation has been used to prove the validity of the mean-field limit for 
a deterministic Cucker-Smale model 
with topological interactions introduced in \cite{H}).

Indeed, given two measures $\rho_1$ and $\rho_2$,
from \eqref{mass} we have
$$
|M_{\rho_1}(B_r(x))-M_{\rho_2}(B_r(x))| \le \| \rho_1 - \rho_2\|_{TV}
$$
where, given $(X, \mathcal{A})$ a measurable space and two measures $\mu$ and $\nu$ over $X$, 
the total variation distance is defined as
$$
\| \mu - \nu\|_{TV}= \sup_{A \in \mathcal{A}} |\mu(A)-\nu(A)|.
$$
In the present work, we use the equivalence between the 
$L^1$ distance and the Total variation for regular measures
and the characterization of the TV distance
given by the Wasserstein distance 
$$
\| \mu - \nu\|_{TV}= \inf_{\pi \in \mathcal{C}(\mu, \nu)} \int_{X \times X} d(x,y)  \de \pi(x,y),
$$
where $\mathcal{C}(\mu, \nu)$ is the set of all couplings, \ie measures on the product space with marginals respectively 
$\mu$ and $\nu$ in the first and second variables,
and $d(a, b) = 1 - \delta_{a,b}$ is the discrete distance (see \cite{V}).

\section{Proof of the result}
\label{sec:proof}

\subsection{Coupling and strategy of the proof}

We introduce, as a coupling between \eqref{generator1} and \eqref{generator2}, the process $ t \to (Z_N(t); \Sigma_N(t))$
on the product space $\R^{2dN} \times \R^{2dN}$, where $\Sigma_N(t)=(Y_N(t), W_N(t))$.
The generator of the new process is  
\[
Q_N=Q_0+\widetilde{Q}_N,
\]
where
\begin{equation}
\label{free-s}
Q_0\Phi(Z_N; \Sigma_N)=
(V_N \cdot \nabla_{X_N} + W_N \cdot \nabla_{Y_N})\Phi(Z_N; \Sigma_N)
\end{equation}
is the free-stream operator, while
\begin{subequations}
\label{coupling}
\begin{align}
\widetilde{Q}_N &\Phi(Z_N; \Sigma_N)= 
\sum_{i=1}^N \sum_{j \neq i} \lambda_{i,j}  [\Phi(X_N, V_N^{i}(v_j); Y_N, W_N^{(i)}(w_j))-\Phi(Z_N; \Sigma_N)] \label{coupling1}\\
&+ \sum_{i=1}^N \sum_{j \neq i}[\pi^N_{i,j}(X_N)- \lambda_{i,j}][\Phi(X_N, V_N^{i}(v_j); \Sigma_N)-\Phi(Z_N; \Sigma_N)] \label{coupling2}\\
&+ \sum_{i=1}^N \sum_{j \neq i}[\pi^{\rho}(y_i,y_j)- \lambda_{i,j}] [\Phi(Z_N; Y_N, W_N^{(i)}(w_j))-\Phi(Z_N; \Sigma_N)] \label{coupling3}\\
&+ \sum_{i=1}^N \int \de u \, \mathcal{E}_i^N(u)[\Phi(Z_N; Y_N, W_N^{(i)}(u))-\Phi(Z_N; \Sigma_N)] \label{coupling4}
\end{align}
\end{subequations}
tends to penalize the discrepancies that can occur over time between $Z_N$ and $\Sigma_N$.

Indeed, in \eqref{coupling1}
the process jumps jointly on both variables 
with a rate given by 
\be
\label{lambda}
\lambda_{i,j}(X_N; y_i, y_j) \coloneq
 \min \{ \pi^N_{i,j}(X_N), \pi^\rho(y_i, y_j)\},
\ee
where
\begin{equation}
\label{massa_rho}
\pi^\rho(y_i, y_j)\coloneq
\alpha_N K\Bigl( M_\rho (B_{|y_i - y_j|}(y_i))\Bigr).
\end{equation}

In \eqref{coupling2} and \eqref{coupling3} the jumps occur only for one of the pair, 
with a transition probability given by the error
between $\lambda_{i,j}$ and $\pi^N$ or $\pi^\rho$.
Finally, in \eqref{coupling4},
\[
\mathcal{E}_i^N(u)=
\int K\Bigl (M_\rho (B_{|y_i - y|}(y_i))\Bigr)f(y, u) \, dy 
- \sum_{j \neq i}\pi^\rho(y_i, y_j) \delta(u - w_j)
\]
is the last error due to the approximation of the limit kinetic equation by the $N$-particle dynamics with transition probabilities given by $\pi^\rho$
and will be treated using the law of large numbers.

We remark that,
since $\int K(x) \de x=1$, formally we have\footnote{In general, the formula is true for $\rho \in L^1(\R^{d})$ and
it is a consequence of the coarea formula
(see \cite[Thm 3.12, p. 140]{EG}).
},
\begin{align*}
\int K\Big(M_\rho(B_{|x-y|}(x)&)\Big) \rho(y) \de y=\int_0^{+\infty} \de r K(M_\rho(B_{r}(x)) \int_{|x-y|=r} \rho(y) \de \mathcal{H}^{n-1}(dy)\\
&= \int_0^{+\infty} \de r K(M_\rho(B_{r}(x)) \frac{\de}{\de r}[M_\rho(B_r(x))]= \int K(x) \de x=1.
\end{align*}
From this fact, it follows that
$Q_N$ is a coupling of the two previously described processes, 
i.e. we recover, 
considering test functions depending only $Z_N$ and $\Sigma_N$ respectively,
the two processes as the two marginals. 

We want to prove that $f$ and $f^N_1$ (defined as in \eqref{marginals})
agree asymptotically in the limit $N \to +\infty$.
To do this we consider $R^N (t)=R^N(Z_N,\Sigma_N;t)$ the law at time $t$ for the coupled process. 
As
initial distribution at time $0$ we assume 
\begin{equation}
\label{initial}
R^N (0)= f_0^{\otimes N}(Z_N) \delta(Z_N - \Sigma_N).
\end{equation}
Let 
$D_N (t)$ 
be the average fraction of particles having different positions or velocities,
\ie using the symmetry of the law,
\begin{equation}
\label{Distance}
D_N(t)=
\int \de R^N(t) \frac{1}{N} \sum_{i=1}^N d(z_i, \sigma_i)= \int \de R^N(t) d(z_1, \sigma_1),
\end{equation}
where $z_i=(x_i, v_i)$, $\sigma_i=(y_i, w_i)$ and 
$d(a, b) = 1 - \delta_{a,b}$ is the discrete distance.

The aim is to show that $D_N(t) \to 0$. 
This means the following: initially the coupled system has all the pairs of particles overlapping.
The dynamics creates discrepancies and the average number of separated pairs 
is exactly $D_N$ which is also the Total Variation distance ($L^1(x,v)$ in our case) between $f^N_1$ and $f$.

Notice that the convergence of the $s$-marginals $f^N_s$ towards $f^{\otimes s}$ claimed in \eqref{main}
is easily recovered by the fact that
\begin{align*}
\|f^{N}_s(t) - f^{\otimes s}(t) \|_{TV} &\le \int \mathfrak{\delta}(Z_s, \Sigma_s) \de R^{N}(Z_N, \Sigma_N; t) \\
& \le \sum_{i=1}^s \int d(z_i, \sigma_i) \de R^N(Z_N, \Sigma_N; t)= sD_N(t)
\end{align*}
where $\delta(a,b)$ denotes the discrete distance on the space $\R^{2ds}\times \R^{2ds}$.

\subsection{Convergence estimates}
Let $S^N_t$ be the semigroup defined by the free-stream generator $Q_0$ in \eqref{free-s}.
To estimate $D_N(t)$ we apply the Duhamel formula in \eqref{Distance} and we get
\begin{multline}
\label{stima1}
 \int \de R^N(t) d(z_1, \sigma_1)=
\int \de R^N(0) d\Bigl(S^N_t(z_1, \sigma_1)\Bigr)\\
+ \int_0^t \de \tau \, \int \de R^N(\tau) \, \widetilde{Q}_N d\Bigl(S^N_{t-\tau}(z_1, \sigma_1)\Bigr),
\end{multline}
where $\widetilde{Q}_N$ is defined in \eqref{coupling}.

The first term in \eqref{stima1} is negligible: 
indeed, from \eqref{initial}, we have
\[
\int \de R^N(0) d\Bigl(S^N_t(z_1, \sigma_1)\Bigr)= 
\int \de f_0^{\otimes N}(Z_N) d\Bigl(S^N_t(z_1, z_1)\Bigr)
\equiv 0.
\]

Concerning the second term in \eqref{stima1}, 
we define 
$$\widebar{z}_1=(x_1+v_1(t-\tau), v_1),\quad  
\widebar{z}_1^{(j)}=(x_1+v_1(t-\tau),v_j)$$ 
and 
$\widebar{X}_N=(x_1+v_1(t-\tau), \dots, x_N+v_N(t-\tau));$
similarly for $\widebar\sigma$, $\widebar{\sigma}^{(j)}$ and $\widebar{Y}_N$.

By \eqref{coupling} we get
\[
\int \de R^N(\tau) \, \widetilde{Q}_N d\Bigl(S^N_{t-\tau}(z_1, \sigma_1)\Bigr)
=A_1(\tau) + A_2(\tau)+A_3(\tau),
\]
where
\[
A_1(\tau)= 
\sum_{j \neq 1} \int \de R^N(\tau) \lambda_{1,j}(\widebar{X}_N;\widebar{y}_1, \widebar{y}_j) 
[d(\widebar{z}_1^{(j)}; \widebar{\sigma}_1^{(j)})-d(\widebar{z}_1; \widebar{\sigma}_1)]
\] 
is due to the term of the generator $\widetilde{Q}_N$ where the velocities of the particles jump simultaneously;
\begin{align*}
A_2(\tau)&=
\sum_{j \neq 1}\int \de R^N(\tau) (\pi^N_{1,j}(\widebar{X}_N) - \lambda_{1,j})
[d(\widebar{z}_1^{(j)}; \widebar{\sigma}_1)-d(\widebar{z}_1; \widebar{\sigma}_1)] \\
&+ \sum_{j \neq 1}\int \de R^N(\tau) (\pi^\rho(\widebar{y}_1,\widebar{y}_j) - \lambda_{1,j})
[d(\widebar{z}_1; \widebar{\sigma}_1^{(j)})-d(\widebar{z}_1; \widebar{\sigma}_1)]
\end{align*}
is due to the terms of the generator where only one of the two coupled processes jump
and
\begin{align*}
A_3(\tau)=
\int \de R^N(\tau) \int \de u \, \widebar{\mathcal{E}}_1^N(u)
[d(\widebar{z}_1; \widebar{\sigma}_1^{(u)})-d(\widebar{z}_1; \widebar{\sigma}_1)]
\end{align*}
is due to the remainder term. Here $\widebar{\mathcal{E}}_1^N(u)$
is $\mathcal{E}_1^N(u)$ evaluated along the moving frame of the free transport.

Here, 
we have used that $d(z_1,\sigma_1)$ depends only on the configurations of the first particle; 
hence, the only non-zero contribution in the sum over $i$ is given for $i=1$.

Concerning $A_1(\tau)$, 
it follows from \eqref{alpha} and \eqref{riemann} that
$$
|e_K(N)| \le \frac{\text{Lip}(K)}{N-1}
$$
and that, for $N>2 \text{Lip(K)}+1$, 
$$
\alpha_N \le \frac{4 \e^{\frac{\text{Lip}(K)}{N-1}}}{N-1},
$$
using the inequality $1/(1-x) \le 4 \e^{x}$ for $x \in (0, 1/2)$.
Therefore, from \eqref{lambda} we get
\[
\lambda_{1,j} \le 
\alpha_N \|K\|_\infty \le \frac{4 \sqrt{\e}\,\text{Lip}(K)}{N-1}.
\]
By the symmetry of $R^N$ and 
denoting $C_K \coloneq 8\sqrt{\e}\,\text{Lip}(K)$, 
\begin{equation}
\label{a1}
A_1(\tau) \le
\frac{C_K}{2(N-1)}  
\sum_{j \neq 1}\int \de R^N(\tau) [d(z_j,\sigma_j)+d(z_1, \sigma_1)]
\le C_K D_N(\tau),
\end{equation}
since 
$d(\widebar{z}_1^{(j)}; \widebar{\sigma}^{(j)}_1)\le d(z_j,\sigma_j)+d(z_1; \sigma_1)$.
Indeed the right-hand side is vanishing 
iff $z_1=\sigma_1$ and $z_j=\sigma_j$ and,
in this case, also the left-hand side is clearly vanishing.

We now give a bound on $A_2(\tau)$. 
Since $\lambda_{1,j}$ 
is the minimum between
$\pi^N_{1,j}$ and $\pi^\rho_{i,j}$, 
we have
\begin{equation}
\label{stima2}
|A_2(\tau)| \le 
\sum_{j \neq 1} \int \de R^N(\tau) 
| \pi^N_{1,j}(\widebar{X}_N) - \pi_{1,j}^\rho(\widebar{y}_1, \widebar{y}_j)|.
\end{equation}
From \eqref{massa_part} and \eqref{massa_rho}, 
\[
| \pi^N_{1,j}(\widebar{X}_N) - \pi_{1,j}^\rho(\widebar{y}_1,\widebar{y}_j)| \le 
\alpha_N \text{Lip}(K)
|M_{\widebar{X}_N}(\widebar{B}^x_{1,j}) - M_\rho(\widebar{B}^y_{1,j})|,
\]
where we are using the shorthand notation
\[
\widebar{B}^x_{1,j}=
B_{|\widebar{x}_1-\widebar{x}_j|}(\widebar{x}_1)
\quad \text{and} \quad 
\widebar{B}^y_{1,j}=B_{|\widebar{y}_1-\widebar{y}_j|}(\widebar{y}_1).
\]
By the triangular inequality
\begin{align*}
|M_{\widebar{X}_N}(\widebar{B}^x_{1,j}) &- M_\rho(\widebar{B}^y_{1,j})|\le 
|M_{\widebar{X}_N}(\widebar{B}^x_{1,j}) - M_{\widebar{X}_N}(\widebar{B}^y_{1,j})|\\
&+|M_{\widebar{X}_N}(\widebar{B}^y_{1,j}) - M_{\widebar{Y}_N}(\widebar{B}^y_{1,j})| 
+|M_{\widebar{Y}_N}(\widebar{B}^y_{1,j}) - M_\rho(\widebar{B}^y_{1,j})|.
\end{align*}
Hence 
we divide the estimate \eqref{stima2} respectively in three terms:
\[
|A_2(\tau)| \le T_1(\tau) + T_2(\tau) + T_3(\tau).
\]

In $T_1(\tau)$ 
we are considering particles with spatial configuration given by $X_N$
and we want to estimate the discrepancy of the configuration
over two different balls $\widebar{B}^x_{1,j}$ and $\widebar{B}^y_{1,j}$.
Since $\widebar{B}^x_{1,j}=\widebar{B}^y_{1,j}$ iff $z_1=\sigma_1$ and $z_j=\sigma_j$,
using that $M_{\widebar{X}_N}\in [0,1]$, we have 
$$
|M_{\widebar{X}_N}(\widebar{B}^x_{1,j}) - M_{\widebar{X}_N}(\widebar{B}^y_{1,j})| \le d(z_1, \sigma_1) + d(z_j, \sigma_j).
$$
Therefore, by the symmetry of $R^N$,
\begin{align*}
T_1(\tau) &\le 
\alpha_N \text{Lip}(K) \sum_{j \neq 1} \int \de R^N(\tau) 
[d(z_1, \sigma_1) + d(z_j, \sigma_j)]\\
&\le C_KD_N(\tau).
\end{align*}

Regarding $T_2(\tau)$, we are considering the discrepancy
of two different configurations over the same ball $\widebar{B}^y_{1,j}$.
Since 
\[
|M_{\widebar{X}_N}(\widebar{B}^y_{1,j}) - M_{\widebar{Y}_N}(\widebar{B}^y_{1,j})|
\le \frac{1}{N} \sum_{i=1}^N d(z_i, \sigma_i),
\]
using again the symmetry of the law, we get
\[
T_2(\tau) 
\le \alpha_N \text{Lip}(K) \sum_{j \neq 1} \int \de R^N(\tau) d(z_1, \sigma_1) 
\le C_K D_N(\tau).
\]

The last estimate on $T_3(\tau)$ is a consequence of the law of large numbers.
After a change of variable, using the symmetry of the law $R^N$
and the fact that
this last term depends only on the $Y_N$ configuration, 
we have that
\[
T_3(\tau) =\alpha_N \text{Lip}(K)
\sum_{j \neq 1} \int \de \rho^{\otimes N}(\tau)|M_{{Y}_N}(B^y_{1,j}) - M_\rho(B^y_{1,j})|,
\]
where $B^y_{1,j}=B_{|y_1-y_j|}(y_1)$.
By Cauchy-Schwartz,
\begin{align*}
\Bigl |&\int \de \rho^{\otimes N}(\tau)|M_{{Y}_N}(B^y_{1,j}) - M_\rho(B^y_{1,j})| \Bigr|^2\\
&\le \int \de \rho^{\otimes N}(\tau)
\Bigl|\frac{1}{N-1} \sum_{h \neq 1} 
\Bigl[\fcr_{B_{1,j}^y}({y}_h) - M_\rho(B_{1,j}^y)\Bigr] \Bigr|^2\\
&\le\sum_{h_1, h_2 \neq 1} \int  \frac{\de\rho^{\otimes N}(\tau)}{(N-1)^2}
\Bigl[\fcr_{B_{1,j}^y}({y}_{h_1}) - M_\rho(B_{1,j}^y)\Bigr]\Bigl [\fcr_{B_{1,j}^y}({y}_{h_2}) - M_\rho(B_{1,j}^y)\Bigr].
\end{align*}
Thanks to the independence of the limit process, 
we get that the only non-zero contributions are given when $h_1=h_2$ and this happens only for $N-1$ terms. 
Hence
\[
T_3(\tau) 
\le \frac{C_K}{\sqrt{N-1}}.
\]
Collecting the estimates on $T_1, T_2$ and $T_3$, 
we obtain that
\begin{equation}
\label{a2}
A_2(\tau) 
\le C_K \Big ( D_N(\tau) + \frac{1}{\sqrt{N-1}}\Big).
\end{equation}

We conclude the proof estimating $A_3(\tau)$. 
Since this term depends only on the independent $Y_N$ configuration
\begin{align*}
|A_3(\tau)|&\le 
\int \frac{\de f^{\otimes N}(\tau)}{N-1}
\sum_{j \neq 1} \Biggl| 
\int K \Bigl(
M_\rho(B_{ |\widebar{y}_1 - y|}(\widebar{y}_1)) 
\Bigr) \de \rho(y) 
- K(M_\rho(\widebar{B}_{1,j}^y)) 
\Biggr|\\
&+ \frac{1}{N-1}\int \de f^{\otimes N}(\tau)  \frac{e_K(N)}{1-e_K(N)}\sum_{j \neq 1}K(M_\rho(\widebar{B}_{1,j}^y)),
\end{align*}
where we added and subtracted the term $\sum_j K(M_\rho(\widebar{B}^y_{1,j}))/(N-1)$.

Applying again the law of large numbers 
on the first term and estimating the second term thanks to
\[
\frac{e_K(N)}{1-e_K(N)} 
\le\frac{C_K}{N-1},
\]
we arrive at
\begin{equation}
\label{a3}
|A_3(\tau)| \le 
\frac{C_K}{\sqrt{N-1}}.
\end{equation}
Collecting the estimates in \eqref{a1}, \eqref{a2} and \eqref{a3} 
and using Gronwall's lemma, we conclude the proof of the theorem.

\section*{acknowledgements} 

PD holds a visiting professor association with the Department of Mathematics, Imperial College London, UK.

\end{document}